\documentclass[12pt]{amsart}
\usepackage{getwidth}
\usepackage{amsfonts}
\usepackage{amsmath}
\usepackage{amssymb}
\usepackage{amsthm}
\newtheorem{theorem}{Theorem}

\newtheorem{lemma}[theorem]{Lemma}
\newenvironment{proof'}{\textbf{Proof}}{\hfill $\square$ \\}

\newcommand{\lf}{L\!f}

\newcommand{\inprod}[2]{\langle #1 \, ,#2 \rangle}
\newcommand{\inprodg}[2]{\langle #1 \, ,#2 \rangle_{G}}
\newcommand{\inprodd}[2]{\langle #1 \, ,#2 \rangle_{D}}
\newcommand{\foral}{\quad \text{for all} \quad}

\newcommand{\sph}{S^{n-1}}

\begin{document}
\numberwithin{theorem}{section} \numberwithin{equation}{theorem}
\pagestyle{plain}
\title{Volumes of Nonnegative Polynomials, Sums of
Squares and Powers of Linear Forms}
\author{Grigoriy Blekherman}
\begin{abstract}
We study the quantitative relationship between the cones of
nonnegative polynomials, cones of sums of squares and cones of
sums of powers of linear forms. We derive bounds on the volumes
(raised to the power reciprocal to the ambient dimension) of
compact sections of the three cones. We show that the bounds are
asymptotically exact if the degree is fixed and number of
variables tends to infinity. When the degree is larger than two it
follows that there are significantly more non-negative polynomials
than sums of squares and there are significantly more sums of
squares than sums of powers of linear forms. Moreover, we quantify
the exact discrepancy between the cones; from our bounds it
follows that the discrepancy grows as the number of variables
increases.
\end{abstract}
\maketitle
\section{Introduction}
Let $P_{n,2k}$ be the vector space of real homogeneous polynomials
in $n$ variables of degree $2k$. There are three interesting
convex cones in $P_{n,2k}$: The cone of nonnegative polynomials,
$C=C_{n,2k}$
\begin{equation*}
C=\bigl{\{}f \in P_{n,2k} \mid f(x) \geq 0 \quad \text{for all}
\quad x \in \mathbb{R}^n \bigr{\}}.
\end{equation*}
The cone of sums of squares, $Sq=Sq_{n,2k}$
\begin{equation*}
Sq=\biggl{\{} f \in P_{n,2k} \mathrel{\bigg{\arrowvert}} f=\sum_i
f_i^2 \quad \text{for some} \quad f_i \in P_{n,k} \biggr{\}}.
\end{equation*}
The cone of sums of $2k$-th powers of linear forms,
$\lf=\lf_{n,2k}$
\begin{equation*}
\lf=\biggl{\{} f \in P_{n,2k} \mathrel{\bigg{\arrowvert}} f=\sum_i
l_i^{2k} \quad \text{for some linear forms} \quad l_i \in
P_{n,1}\bigg{\}}.
\end{equation*}
A different notation of $P_{n,2k},\Sigma_{n,2k}$ and $Q_{n,2k}$
respectively was employed by Reznick in the study of these cones
\cite{rez}. The cones are clearly nested:
\begin{equation*}
\lf_{n,2k} \subseteq Sq_{n,2k} \subseteq C_{n,2k}.
\end{equation*}
\indent It is known that for quadratic forms these cones coincide.
Moreover, it is not hard to show that in all other cases there are
sums of squares that are not $2k$-th powers of linear forms.
Hilbert proved that in the cases $n=2$, $k=1$ and, $n=3$ and
$k=2$, a nonnegative polynomial is necessarily a sum of squares;
in all other cases there exist nonnegative polynomials that are
not sums of squares \cite{hilbert}. The situation with respect to
containment has therefore been
completely known for a long time.\\
\indent There remains, however, the question of the quantitative
relationship between these cones. There are several known families
of polynomials that are not sums of squares \cite{choi},
\cite{rez4}; however all of these examples lie close to the
boundary of the cone of nonnegative polynomials. To the author's
knowledge little except for the equality in the case of quadratic
forms is known. In this paper we show that the picture is quite
different for a fixed degree greater than 2.
\\ \indent For a convex set $K$ a good measure of size of $K$ that
takes into account the effect of large dimensions is the volume of
$K$ raised to the power reciprocal to the ambient dimension:
\begin{equation*}
(\text{Vol} \, K)^{1/\text{dim} \, K}.
\end{equation*}
For example, homothetically expanding $K$ by a constant factor
leads to an increase by the same factor in this normed volume.
\\ \indent We derive bounds on volumes, raised to the power reciprocal to
the ambient dimension, of sections of the three cones with the
hyperplane of all forms of integral 1 on the unit sphere $S^{n-1}$
in $\mathbb{R}^n$. We show that the bounds are asymptotically
tight if the degree is fixed and number of variables tends to
infinity. If the degree is greater than 2 then the order of
dependence on the number of variables $n$ is quite different for
the three cones. We remark that this indeed shows that
asymptotically the cones differ drastically in size. These bounds
provide us with the complete picture of metric dependence of the
size of all three cones on the number of variables, when the
degree is fixed.
\\ \indent
We would also like to mention that the bounds that separate the
cone of nonnegative polynomials from the cone of sums of squares
are interesting from the point of view of computational complexity
\cite{parrilo}. Namely, they show that it is not feasible in
general to replace testing for positivity with testing whether a
polynomial is a sum of squares, since for degree greater than two
the sizes of the cones are drastically different. Some of the
bounds given in this paper have already been proved by the author
in \cite{greg}; we reproduce their proofs for the sake of
completeness.
\section{Main Theorems}

We begin by introducing some notation. In order to compare the
cones we take compact bases. Let $M=M_{n,2k}$ be the hyperplane of
all forms in $P_{n,2k}$ with integral 0 on the unit sphere $\sph$:
\begin{equation*}
M_{n,2k}=\biggr{\{} f \in P_{n,2k} \mathrel{\bigg{\arrowvert}}
\int_{\sph}f \, d\sigma=0 \biggl{\}}.
\end{equation*}
Let $r^{2k}$ in $P_{n,2k}$ be the polynomial constant on the unit
sphere $\sph$:
\begin{equation*}
r^{2k}=(x_1^2+ \ldots +x_n^2)^k.
\end{equation*}
\indent Let $M'$ be the affine hyperplane of all forms of integral
1 on the unit sphere $\sph$. We define compact convex bodies
$\widetilde{C}$, $\widetilde{Sq}$ and $\widetilde{\lf}$ by
intersecting the respective cones with $M'$ and then translating
the compact intersection into $M$ by subtracting $r^{2k}$.
Formally we can define $\widetilde{C}$, $\widetilde{Sq}$ and
$\widetilde{\lf}$ as the sets of all forms $f$ in $M_{n,2k}$ such
that $f+r^{2k}$ lies in the respective cone:
\begin{eqnarray*}
\widetilde{C}=\{f \in M_{n,2k} \quad \mid \quad f+r^{2k} \in C \}, \\
\widetilde{Sq}=\{f \in M_{n,2k} \quad \mid \quad f+r^{2k} \in Sq \}, \\
\widetilde{\lf}=\{f \in M_{n,2k} \quad \mid \quad f+r^{2k} \in \lf
\}.
\end{eqnarray*}
We note that these sections are the natural ones to take since
$M_{n,2k}$ is the only
linear hyperplane in $P_{n,2k}$ that is preserved by an orthogonal change of coordinates in $\mathbb{R}^n$.\\
\indent We work with the following Euclidean metric on $P_{n,2k}$,
which we call the integral or $L^2$ metric,
\begin{equation*}
\inprod{f}{g}=\int_{\sph} fg \, d\sigma,
\end{equation*}
where $\sigma$ is the rotation invariant probability measure on
$\sph$. We use $D_M$ to denote the dimension of $M_{n,2k}$, $S_M$
to denote the unit sphere in $M_{n,2k}$ and $B_M$ to denote the
unit ball in $M_{n,2k}$. The main results of this paper are the
following three theorems:
\begin{theorem}
\label{posmain} There exist constants $\alpha_1$ and $\beta_1$$>0$
dependent only on $k$ such that
\begin{equation*}
\beta_1 n^{-1/2} \leq \bigg{(}\frac{\text{Vol} \,
\widetilde{C}}{\text{Vol} \, B_M}\bigg{)}^{1/D_M} \leq
\alpha_1n^{-1/2}.
\end{equation*}
\end{theorem}
\begin{theorem}
\label{squarenormbound} There exist constants $\alpha_2$ and
$\beta_2$$>0$ dependent only on $k$ such that
\begin{equation*}
\beta_2n^{-k/2} \bigg{(}\frac{\text{Vol}\,
\widetilde{Sq}}{\text{Vol}\, B_M}\bigg{)}^{1/D_M} \leq \alpha_2
n^{-k/2}.
\end{equation*}
\end{theorem}

\begin{theorem}
\label{powersnormbound} There exist constants $\alpha_3$ and
$\beta_3$ $>0$ dependent only on $k$ such that for all $\epsilon
> 0$ and $n$ large enough
\begin{equation*}
\beta_3n^{-k+1/2} \leq \bigg{(}\frac{\text{Vol} \,
\widetilde{Lf}}{\text{Vol} \, B_M}\bigg{)}^{1/D_M}   \leq
\alpha_3n^{-k+1/2+\epsilon}.
\end{equation*}
\end{theorem}
\indent We observe that if the degree $2k$ is equal to two, then
all of the above bounds agree asymptotically. However if the
degree is greater than two then we see that the bases
$\widetilde{C}$, $\widetilde{Sq}$ and $\widetilde{\lf}$
asymptotically have quite
different volumes.\\
\indent The rest of the paper is structured as follows. In Section
3 we collect preliminary material necessary for the proofs. Since
many of the estimates used are technical in nature, in Section 4
we give an outline of the proofs postponing the technical details
for the later sections. In Section 5 we prove the bounds for the
cone of nonnegative polynomials. In Section 6 we introduce a
different metric on $P_{n,2k}$ and prove duality results used
later on. In Section 7 we prove the bounds for the cone of sums of
squares and in Section 8 we prove the bounds for the cone of sums
of powers of linear forms.
\section{Preliminaries}
\subsection{The Action of the Orthogonal Group on $P_{n,2k}$}\

There is the following action of $SO(n)$ on $P_{n,2k}$,
\begin{equation*}
A \in SO(n) \quad \text{sends} \quad f \in  P_{n,2k} \quad
\text{to} \quad Af=f(A^{-1}x).
\end{equation*}
We observe that the cones $C$, $Sq$ and $Lf$ are invariant under
this action and so is $M_{n,2k}$, the hyperplane of polynomials of
integral $0$. Therefore the sections $\widetilde{C}$,
$\widetilde{Sq}$ and $\widetilde{\lf}$ are fixed by $SO(n)$ as
well.
\\ \indent Let $\Delta$ be the Laplace differential
operator:
\begin{equation*}
\Delta=\frac{\partial^2}{\partial x_1^2}+\ldots
+\frac{\partial^2}{\partial x_n^2}.
\end{equation*}
A form $f$ such that
\begin{equation*}
\Delta(f)=0,
\end{equation*}
is called \textit{harmonic}. We will need the fact that the
irreducible components of this representation are subspaces
$H_{n,2l}$ for $0\leq l \leq k$, which have the following form:
\begin{equation*}
H_{n,2l}=\bigl{\{} f \in P_{n,2k} \mid f=r^{2k-2l}h \quad
\text{where} \quad h \in P_{n,2l} \quad \text{is harmonic}
\bigr{\}}.
\end{equation*}
\indent For $v \in \mathbb{R}^n$ the functional
\begin{equation*}
\lambda_v:M_{n,2k} \longrightarrow \mathbb{R}, \qquad
\lambda_v(f)=f(v),
\end{equation*}
is linear and therefore there exists a form $q_v \in M$ such that
\begin{equation*}
\lambda_v(f)=\inprod{q_v}{f}.
\end{equation*}
There are explicit descriptions of the polynomials $q_v$, under a
suitable normalization they are so called Gegenbauer or
ultraspherical polynomials. We will only need the property that
for $v \in \sph$
\begin{equation*}
||\,q_v||_{2}=\sqrt{D_M}.
\end{equation*}
For more details on this representation of $SO(n)$ see
\cite{vilenkin}.

\subsection{The Blaschke-Santal\'{o} Inequality}\hspace{.1cm}
\\\indent
Let $K$ be a full-dimensional convex body in $\mathbb{R}^n$ with
origin in its interior and let $\inprod{\,}{\,}$ be an inner
product. We will use $K^{\circ}$ to denote the polar of $K$,
\begin{equation*}
K^{\circ}=\bigl{\{} x \in \mathbb{R}^n \, \mid \, \inprod{x}{y}
\leq 1 \foral y \in K \bigr{\}}.
\end{equation*}
Now suppose that a point $z$ is in the interior of $K$ and let
$K^z$ be the polar of $K$ when $z$ is translated to the origin:
\begin{equation*}
K^z=\bigl{\{} x \in \mathbb{R}^n \, \mid \, \inprod{x-z}{y-z} \leq
1 \foral y \in K \bigr{\}}.
\end{equation*}
The point $z$ at which the volume of $K^z$ is minimal is unique
and it is called the Santal\'{o} point of $K$. Moreover the
following inequality on volumes of $K$ and $K^z$ holds:
\begin{equation*}
\frac{\text{Vol} \, K \, \text{Vol} \, K^z}{(\text{Vol} \, B)^2}
\leq 1,
\end{equation*}
where $B$ is the unit ball of $\inprod{\,}{\,}$ and $z$ is the
Santal\'{o} point of $K$. This is known as the
Blaschke-Santal\'{o} inequality \cite{santalo}.

\section{Outline of Proofs}
Since many of the following proofs are technical we would like to
first give an informal outline.
\\ \indent We begin with the description of the proofs for the cone of
nonnegative polynomials. We observe that $\widetilde{C}$ is the
convex body of forms of integral $0$ on $\sph$, such that the
minimum of the forms on $\sph$ is at least $-1$,
\begin{equation*}
\widetilde{C}=\big{\{}f\in M_{n,2k} \quad \mid \quad f(x) \geq -1
\foral x \in\sph \big{\}}.
\end{equation*}
Let $B_{\infty}$ be the unit ball of $L^{\infty}$ norm in
$M_{n,2k}$,
\begin{equation*}
B_{\infty}=\big{\{}f\in M_{n,2k} \quad \mid \quad |f(x)| \leq 1
\foral x \in \sph \big{\}}.
\end{equation*}
It follows that
\begin{equation*}
B_{\infty}=\widetilde{C} \cap -\widetilde{C} \qquad \text{and
therefore} \qquad B_{\infty} \subset \widetilde{C}.
\end{equation*}
However, using the Blaschke-Santal\'{o} inequality and a theorem
of\\ Rogers and Shephard \cite{pach} we can show that conversely
\begin{equation*}
\bigg{(}\frac{\text{Vol} \, B_{\infty}}{\text{Vol} \,
\widetilde{C}}\bigg{)}^{1/D_M} \geq 1/4.
\end{equation*}
Therefore it suffices to derive upper and lower bounds for the
volume of $B_{\infty}$.
\\ \indent For the lower bound we reduce the proof to bounding
the average $L^{\infty}$ norm of a form in $M_{n,2k}$,
\begin{equation*}
\int_{S_M} ||f||_{\infty}\, d\mu,
\end{equation*}
where $S_M$ is the unit sphere in $M_{n,2k}$ and $\mu$ is the
rotation invariant probability measure on $S_M$. The key idea is
to estimate $||f||_{\infty}$ using $L^{2p}$ norms for some large
$p$. An inequality of Barvinok \cite{barv} is used to see that
taking $p=n$ suffices for $||f||_{2p}$ to be within a constant
factor of $||f||_{\infty}$. The proof is completed with some
estimates.
\\ \indent The techniques used for the proof of the upper bound
are quite different. Let $\nabla f$ be the gradient of $f \in
P_{n,2k}$,
\begin{equation*}
\nabla f = \bigg{(} \frac{\partial f}{\partial x_1}\, ,\ldots , \,
\frac{\partial f}{\partial x_n} \bigg{)},
\end{equation*}
and let $\inprod{\nabla{f}}{\nabla{f}}$ be the following
polynomial giving the squared length of the gradient of $f$,
\begin{equation*}
\inprod{\nabla f}{\nabla f}= \bigg{(}\frac{\partial f}{\partial
x_1}\bigg{)}^2  + \ldots + \bigg{(}\frac{\partial f}{\partial x_n}
\bigg{)}^2.
\end{equation*}
\indent The key to the proof is the following theorem of Kellogg
\cite{kellogg} which tells us that for homogeneous polynomials the
maximum length of the gradient on the unit sphere $\sph$ is equal
to the maximum absolute value of the polynomial on $\sph$
multiplied by the degree of the polynomial:
\begin{equation*}
||\inprod{\nabla f}{\nabla f}||_{\infty}=4k^2||f||^2_{\infty}.
\end{equation*}
Now we define a different inner product on $P_{n,2k}$ which we
call the gradient inner product,
\begin{equation*}
\inprodg{f}{g}=\frac{1}{4k^2}\int_{\sph} \inprod{\nabla f}{\nabla
g} \, d\sigma.
\end{equation*}
We denote the norm of $f$ in the gradient metric by $||f||_G$ and
the unit ball of the gradient metric in $M_{n,2k}$ by $B_G$. We
observe that
\begin{equation*}
||f||_G=\frac{1}{4k^2}\int_{\sph} \inprod{\nabla f}{\nabla f} \,
d\sigma,
\end{equation*}
and hence it follows that
\begin{equation*}
||f||_G \leq ||f||_{\infty} \qquad \text{and therefore} \qquad
B_{\infty} \subset B_G.
\end{equation*}
\indent The relationship between the gradient metric and the
integral metric can be calculated precisely by using the fact that
both metrics are $SO(n)$-invariant. Therefore these metrics are
constant multiples of each other in the irreducible subspaces of
the $SO(n)$ representation and the constants can be calculated
directly using the Stokes' formula. Hence we obtain an upper bound
for the volume of $B_{\infty}$ in terms of the volume of $B_M$,
the unit ball of the $L^2$ metric in $M_{n,2k}$.
\\ \indent The intuitive idea of the proof is as follows. In the
$L^2$ metric we have,
\begin{equation*}
||f||_2 \leq ||f||_{\infty} \qquad \text{and therefore} \qquad
B_{\infty} \subset B_M.
\end{equation*}
However we give up too much in this estimate. On the other hand,
it is not hard to show that
\begin{equation*}
f^2(x) \leq 4k^2 \inprod{\nabla f}{\nabla f} \foral x \in \sph.
\end{equation*}
Direct computations show that using the gradient metric gives us a
better estimate and that this estimate is fine enough for our
purposes.
\\ \indent The proof of the upper bound for the cone of sums of
squares is quite similar to the proof of the lower bound for the
cone of nonnegative polynomials. We define the following norm on
$P_{n,2k}$,
\begin{equation*}
||f||_{sq}=\max_{g \in S_{P_{n,k}}} |\inprod{f}{g^2}|,
\end{equation*}
where $S_{P_{n,k}}$ is the unit sphere in $P_{n,k}$. Using
inequalities from convexity we can reduce the proof to bounding
the average $||f||_{sq}$.
\\ \indent To every form $f \in P_{n,2k}$ we can associate a
quadratic form $H_f$ on $P_{n,2k}$ by letting
\begin{equation*}
H_f(g)=\inprod{f}{g^2} \qquad \text{for} \qquad g \in P_{n,k}.
\end{equation*}
It follows that
\begin{equation*}
||f||_{sq}=||H_f||_{\infty}.
\end{equation*}
Now we can estimate $||H_f||_{\infty}$ by high $L^{2p}$ norms of
$H_f$ and the proof is finished using similar ideas to the proof
for the case of nonnegative polynomials.
\\ \indent For the remainder of the proofs we will need to consider
yet another metric on $P_{n,2k}$. To a form $f \in P_{n,2k}$,
\begin{equation*}
f=\sum_{\alpha=(i_1, \ldots ,i_n)}c_{\alpha}x_1^{i_1}\ldots
x_n^{i_n}.
\end{equation*}
we formally associate the differential operator $D_f$:
\begin{equation*}
D_f=\sum_{\alpha=(i_1, \ldots
,i_n)}c_{\alpha}\frac{\partial^{i_1}}{\partial x_1^{i_1}}\cdots
\frac{\partial^{i_n}}{\partial x_n^{i_n}}.
\end{equation*}
We define the following metric on $P_{n,2k}$, which we call the
differential metric:
\begin{equation*}
\inprodd{f}{g}=D_f(g).
\end{equation*}
It is not hard to check that this indeed defines a symmetric
positive definite bilinear form, which is invariant under the
action of $SO(n)$. The relationship between the differential
metric and the integral metric can be calculated precisely.
\\ \indent For the proof of the lower bound for the cone of sums of
squares we show that the dual cone $Sq^*_d$ of $Sq$ with respect
to the differential metric is contained in $Sq$. Therefore we can
derive a lower bound on the volume of $\widetilde{Sq}$ by using
the Blaschke-Santal\'{o} inequality.
\\ \indent It can be shown that the cone of sums of $2k$-th powers
of linear forms $Lf$ is dual to $C$ in the differential metric.
The proofs of the bounds follow from the bounds derived for
$\widetilde{C}$ and the Blaschke-Santal\'{o} inequality.

\section{Nonnegative Polynomials}
\setcounter{equation}{0} In this section we prove Theorem
\ref{posmain}. Here is the precise statement of the bounds:
\begin{theorem}
\label{posmainfull} There are the following bounds on the volume
of $\widetilde{C}$:
\begin{equation*}
\frac{1}{2\sqrt{4k+2}} \, n^{-1/2} \leq \bigg{(}\frac{\text{Vol}
\, \widetilde{C}}{\text{Vol} \, B_M}\bigg{)}^{1/D_M} \leq
4\bigg{(}\frac{2k^2}{4k^2+n-2}\bigg{)}^{1/2}.
\end{equation*}
\end{theorem}
\subsection{Proof of the Lower Bound} \hspace{1mm}\\
\indent For a real Euclidean vector space $V$ with the unit sphere
$S_V$ and a function $f:V \rightarrow \mathbb{R}$ we use
$||f||_{p}$ to denote the $L^p$ norm of $f$:
\begin{equation*}
||f||_p=\bigg{(}\int_{S_V} |f|^p \, d\mu \bigg{)}^{1/p} \qquad
\text{and} \qquad ||f||_{\infty}=\max_{x \in S_V} |f(x)|.
\end{equation*}
\indent We begin by observing that $\widetilde{C}$ is a convex
body in $M_{n,2k}$ with origin in its interior and the boundary of
$\widetilde{C}$ consists of polynomials with minimum $-1$ on
$\sph$. Therefore the gauge $G_C$ of $\widetilde{C}$ is given by:
\begin{equation*}
G_{C}(f)=|\min_{v \, \in \sph} f(v)\,|.
\end{equation*}
By using integration in polar coordinates in $M$ we obtain the
following expression for the volume of $\widetilde{C}$,
\begin{equation}
\label{integralvolume} \biggl{(}\frac{\text{Vol} \,
\widetilde{C}}{\text{Vol }
B_M}\biggr{)}^{\frac{1}{D_M}}=\biggl{(}\int_{S_M} G_{C}^{-D_M} \,
d\mu \biggr{)}^{\frac{1}{D_M}},
\end{equation}
where $\mu$ is the rotation invariant probability measure on
$S_M$. The relationship \eqref{integralvolume}
holds for any convex body with origin in its interior \cite[p. 91]{pisier}. \\
We interpret the right hand side of \eqref{integralvolume} as
$||G_C^{-1}||_{D_M}$, and by H\"{o}lder's inequality
\begin{equation*}
||G_C^{-1}||_{D_M} \geq ||G_C^{-1}||_1.
\end{equation*}
Thus,
\begin{equation*}
\biggl{(}\frac{\text{Vol} \, \widetilde{C}}{\text{Vol} \,
B_M}\biggr{)}^{\frac{1}{D_M}} \geq \int_{S_M} G_C^{-1} \, d\mu.
\end{equation*}
By applying Jensen's inequality \cite[p.150]{hardy}, with convex
function $y=1/x$ it follows that,
\begin{equation*}
\int_{S_M} G_C^{-1} \ d\mu \geq \bigg{(}\int_{S_M} G_C \, d\mu
\bigg{)}^{-1}.
\end{equation*}
Hence we see that
\begin{equation*}
\bigg{(} \frac{\text{Vol} \, \widetilde{C}}{\text{Vol} \, B_M}
\bigg{)}^{\frac{1}{D_M}} \geq \bigg{(}\int_{S_M} |\min f| \, d \mu
\bigg{)}^{-1}.
\end{equation*}
Clearly, for all $f \in P_{n,2k}$
\begin{equation*}
||f||_{\infty} \geq |\min f|.
\end{equation*}
Therefore,
\begin{equation*}
\bigg{(} \frac{\text{Vol} \, \widetilde{C}}{\text{Vol} \, B_M}
\bigg{)}^{\frac{1}{D_M}} \geq \bigg{(}\int_{S_M} ||f||_{\infty} \,
d \mu \bigg{)}^{-1}.
\end{equation*}
The proof of the lower bound of Theorem \ref{posmainfull} is now
completed by the following estimate.
\begin{theorem}
\label{infinitynorm} Let $S_M$ be the unit sphere in $M_{n,2k}$
and let $\mu$ be the rotation invariant probability measure on
$S_M$. Then the following inequality for the average $L^{\infty}$
norm over $S_M$ holds:
\begin{equation*}
\int_{S_M} ||f||_{\infty} \, d\mu \leq 2\sqrt{2n(2k+1)}.
\end{equation*}
\end{theorem}
\begin{proof}
It was shown by Barvinok in \cite{barv} that for all $f \in
P_{n,2k}$,
\begin{equation*}
||f||_{\infty} \leq \binom{2kn+n-1}{2kn}^{\frac{1}{2n}}||f||_{2n}.
\end{equation*}
By applying Stirling's formula we can easily obtain the bound
\begin{equation*}
\binom{2kn+n-1}{2kn}^{\frac{1}{2n}} \leq 2\sqrt{2k+1}.
\end{equation*}
Therefore it suffices to estimate the average $L^{2n}$ norm, which
we denote by $A$:
\begin{equation*}
A=\int_{S_M} ||f||_{2n} \, d\mu.
\end{equation*}
Applying H\"{o}lder's inequality we observe that
\begin{equation*}
A=\int_{S_M} \bigg{(}\int_{\sph} f^{2n}(x) \, d\sigma
\bigg{)}^{\frac{1}{2n}} d\mu \leq \bigg{(}\int_{S_M}\int_{\sph}
f^{2n}(x)\, d\sigma \, d\mu \bigg{)}^{\frac{1}{2n}}.
\end{equation*}
By interchanging the order of integration we obtain
\begin{equation}
\label{first} A \leq \bigg{(} \int_{\sph} \int_{S_M} f^{2n}(x) \,
d\mu \, d\sigma \bigg{)}^{\frac{1}{2n}}.
\end{equation}
We now note that by symmetry of $M$
\begin{equation*}
\int_{S_M} f^{2n}(x) \, d\mu,
\end{equation*}
is the same for all $x \in \sph$. Therefore we see that in
\eqref{first} the outer integral is redundant and thus
\begin{equation}
\label{second} A \leq \bigg{(} \int_{S_M} f^{2n}(v) \, d\mu
\bigg{)}^{\frac{1}{2n}}, \qquad \text{where} \ v \ \text{is any
vector in} \ \sph.
\end{equation}
\indent We recall from Section 3 that for $v \in \sph$ there there
exists a form $q_v$ in $M$ such that
\begin{equation*}
\inprod{f}{q_v}=f(v) \foral f \in M \qquad \text{and} \qquad
||q_v||_2=\sqrt{D_M}.
\end{equation*}

Rewriting \eqref{second} we see that
\begin{equation}
\label{third} A \leq \bigg{(} \int_{S_M} \inprod{f}{q_v}^{2n} \,
d\mu \bigg{)}^{\frac{1}{2n}}.
\end{equation}
We observe that
\begin{equation*}
\int_{S_M} \inprod{f}{q_v}^{2n} \, d\mu=(D_M)^n \
\frac{\Gamma(n+\frac{1}{2}\,) \,
\Gamma(\frac{1}{2}D_M)}{\sqrt{\pi} \, \Gamma(\frac{1}{2}D_M+n)}.
\end{equation*}
We substitute this into \eqref{third} to obtain,
\begin{equation*}
A \leq \bigg{(}(D_M)^n \ \frac{\Gamma(n+\frac{1}{2}\,) \,
\Gamma(\frac{1}{2}D_M)}{\sqrt{\pi} \, \Gamma(\frac{1}{2}D_M+n)}
\bigg{)}^{\frac{1}{2n}}.
\end{equation*}
Since
\begin{equation*}
\bigg{(}\frac{\Gamma(\frac{1}{2}D_M)}{\Gamma(\,\frac{ 1}{2}D_M
+n)}\bigg{)}^{\frac{1}{2n}} \leq \sqrt{\frac{2}{D_M}} \qquad
\text{and} \qquad \bigg{(}
\frac{\Gamma(n+1/2\,)}{\sqrt{\pi}}\bigg{)}^{\frac{1}{2n}} \leq
n^{1/2},
\end{equation*}
we see that
\begin{equation*}
A \leq (2n)^{1/2}.
\end{equation*}
The theorem now follows.
\end{proof}
\subsection{Proof of the Upper Bound} \hspace{1mm}\\
\indent We begin by noting that the origin is the only point in
$M$ fixed by $SO(n)$. Let $\widetilde{C}^{\circ}$ be the polar of
$\widetilde{C}$ in $M_{n,2k}$,
\begin{equation*}
\widetilde{C}^{\circ}=\{f \in M_{n,2k} \mid \inprod{f}{g} \leq 1
\foral g \in \widetilde{C}\}.
\end{equation*}
Since $\widetilde{C}$ is fixed by the action of $SO(n)$ and
Santal\'{o} point of a convex body is unique, it follows that the
origin is the Santal\'{o} point of $\widetilde{C}$. We now use
Blaschke-Santal\'{o} inequality, which applied to $\widetilde{C}$
gives us:
\begin{equation*}
(\text{Vol} \, \widetilde{C}) \, (\text{Vol}\,
\widetilde{C}^{\circ})\leq (\text{Vol}\, B_M)^2.
\end{equation*}
Therefore it would suffice to show that
\begin{equation}
\label{ratio} \bigg{(}\frac{\text{Vol} \,
\widetilde{C}^{\circ}}{\text{Vol} \, B_M}\bigg{)}^{1/D_M} \geq
\frac{1}{4}\bigg{(}\frac{4k^2+n-2}{2k^2}\bigg{)}^{1/2}.
\end{equation}

\indent Let $B_{\infty}$ be the unit ball of the $L^{\infty}$
metric in $M_{n,2k}$,
\begin{equation*}
B_{\infty}=\{f \in M \ \mid \ ||f||_{\infty} \leq 1 \}.
\end{equation*}
We observe that $B_{\infty}$ is clearly the intersection of
$\widetilde{C}$ with $-\widetilde{C}$:
\begin{equation*}
B_{\infty}=\widetilde{C} \cap -\widetilde{C}.
\end{equation*}
By taking polars it follows that
\begin{equation*}
B_{\infty}^{\circ} =
\text{ConvexHull}\{C^{\circ},-\,C^{\circ}\}\subset \,
\widetilde{C}^{\circ} \oplus(-\, \widetilde{C}^{\circ}),
\end{equation*}
where $\, \oplus \, $ denotes Minkowski addition. By theorem of
Rogers and Shephard, \cite{pach} p. 78, it follows that
\begin{equation*}
\text{Vol} \, B_{\infty}^{\circ} \leq \binom{2D_M}{D_M} \text{Vol}
\, \widetilde{C}^{\circ}.
\end{equation*}
Since
\begin{equation*}
\binom{2D_M}{D_M} \leq 4^{D_M},
\end{equation*}
we obtain
\begin{equation*}
\bigg{(}\frac{\text{Vol} \, \widetilde{C}^{\circ}}{\text{Vol} \,
B_{\infty}^{\circ}}\bigg{)}^{1/D_M} \geq \frac{1}{4}.
\end{equation*}
Combining with \eqref{ratio} we see that we have reduced the lower
bound of Theorem \ref{posmainfull} to showing that
\begin{equation}
\label{reduce} \bigg{(}\frac{\text{Vol} \,
B_{\infty}^{\circ}}{\text{Vol} \, B_M}\bigg{)}^{1/D_M} \geq
\bigg{(}\frac{4k^2+n-2}{2k^2}\bigg{)}^{1/2}
\end{equation}

\indent For a form $f$ we use $\nabla f$ to denote the gradient of
$f$:
\begin{equation*}
\nabla f =\bigg{(} \frac{\partial f}{\partial x_1}\, ,\ldots \, ,
\frac{\partial f}{\partial x_n} \bigg{)}.
\end{equation*}
We also define a different Euclidean metric on $P_{n,2k}$ which we
call the gradient metric:
\begin{equation*}
\inprodg{f}{g}=\frac{1}{4k^2}\int_{\sph} \inprod{\nabla f}{\nabla
g} \, d\sigma.
\end{equation*}
We denote the unit ball in this metric by $B_G$ and the norm of
$f$ by $||f||_G$. For $f \in P_{n,2k}$ let $\inprod{\nabla
f}{\nabla f}$ be the following polynomial:
\begin{equation*}
\inprod{\nabla f}{\nabla f}= \bigg{(}\frac{\partial f}{\partial
x_1}\bigg{)}^2  + \ldots + \bigg{(}\frac{\partial f}{\partial x_n}
\bigg{)}^2.
\end{equation*}
It was shown by Kellogg in \cite{kellogg} that
\begin{equation*}
||\inprod{\nabla f}{\nabla f}||_{\infty}=4k^2||f||^2_{\infty}.
\end{equation*}
It clearly follows that
\begin{equation*}
||f||_{\infty} \geq ||f||_G,
\end{equation*}
and therefore
\begin{equation*}
B_{\infty} \subseteq B_G.
\end{equation*}
Polarity reverses inclusion and thus we see that
\begin{equation*}
B_{G}^{\circ} \subseteq B_{\infty}^{\circ} \quad \text{and} \quad
\text{Vol} \, B_G^{\circ}=\frac{(\text{Vol} \, B_M)^2}{\text{Vol}
\, B_G},
\end{equation*}
since $B_G$ is an ellipsoid. Thus \eqref{reduce} and consequently
the upper bound of Theorem \ref{posmainfull} will follow from the
following lemma.
\begin{lemma}
\begin{equation*}
\bigg{(} \frac{\text{Vol} \, B_{M}}{\text{Vol} \, B_{\,G}}
\bigg{)}^{1/D_M} \geq \bigg{(}\frac{4k^2+n-2}{2k^2}\bigg{)}^{1/2}.
\end{equation*}
\end{lemma}
\begin{proof}
It will suffice to show that for all $f \in M$
\begin{equation}
\label{product} \inprodg{f}{f} \geq
\frac{4k^2+n-2}{2k^2}\inprod{f}{f}.
\end{equation}
\indent By the invariance of both inner products under the action
of $SO(n)$, it is enough to prove \eqref{product} in the
irreducible
components of the representation.\\
\indent First let $f$ be a harmonic form of degree $2d$ in $n$
variables. Then we claim that
\begin{equation*}
\inprod{f}{f}=\frac{2d}{4d+n-2}\, \inprodg{f}{f}.
\end{equation*}
Indeed consider the vector field $F=f(v)\, \! \nabla \! f$ on
$\sph$. By the Divergence Theorem:
\begin{equation*}
\int_{\sph} \inprod{F}{v} \, dx(v) = \int_{||x|| \leq 1}
\text{div}\, F \, dx,
\end{equation*}
where $dx$ is the Lebesgue measure and $\text{div} \, F$ is the
divergence of $F$:
\begin{equation*}
\text{div} \, F = \frac{\partial{F_1}}{\partial{x_1}} + \ldots +
\frac{\partial{F_n}}{\partial{x_n}}.
\end{equation*}
Since $f$ is homogeneous of degree $2d$, it follows that
\begin{equation*}
\inprod{\nabla f}{v}=2d\,f(v).
\end{equation*}
Therefore
\begin{equation*}
\int_{\sph} \inprod{F}{v} \, dx=2\omega_n d\int_{\sph} f^2 \,
d\sigma=2\omega_n d \inprod{f}{f},
\end{equation*}
where $\omega_n$ is the surface area of $\sph$. Since $f$ is
harmonic it follows that
\begin{equation*}
\text{div} \,
F=\bigg{(}\frac{\partial{f}}{\partial{x_1}}\bigg{)}^2 + \ldots +
\bigg{(}\frac{\partial{f}}{\partial{x_n}}\bigg{)}^2=\inprod{\nabla
f}{\nabla f}.
\end{equation*}
We observe that $\inprod{\nabla f}{\nabla f}$ is a homogeneous
polynomial of degree $4d-2$ and therefore
\begin{equation*}
\int_{||x|| \leq 1} \inprod{\nabla f}{\nabla f} \,
dx=\frac{\omega_n}{4d+n-2}\int_{\sph} \inprod{\nabla f}{\nabla
f}\, d\sigma.
\end{equation*}
The claim now follows. \\
\indent Now suppose that $f=hr^{2k-2d}$ where $h$ is a harmonic
form of degree $2d \leq 2k$. It is easy to check that
\begin{equation*}
\inprodg{f}{f}=\frac{d^2}{k^2}\inprodg{h}{h}+\frac{k^2-d^2}{k^2}\inprod{h}{h}.
\end{equation*}
We know that
\begin{equation*}
\inprodg{h}{h}=\frac{4d+n-2}{2d}\inprod{h}{h} \quad \text{and}
\quad \inprod{f}{f}=\inprod{h}{h}.
\end{equation*}
Thus
\begin{equation*}
\inprodg{f}{f}=\frac{2k^2+d(n-2)+2d^2}{2k^2}\inprod{f}{f}.
\end{equation*}
Since $f \in M_{n,2k}$ we know that $1\leq d \leq k$. The minimum
clearly occurs when $d=1$ and we see that
\begin{equation*}
\inprodg{f}{f} \leq \frac{4k^2+n-2}{2k^2}\inprod{f}{f}.
\end{equation*}
The lemma now follows.
\end{proof}
\section{The Differential Metric}

\setcounter{equation}{0} Before we proceed with the proofs of
Theorems \ref{squarenormbound} and \ref{powersnormbound} we will
need some preparatory results that involve switching to a
different Euclidean metric on $P_{n,2k}$. \\
\indent To a form $f \in P_{n,2k}$,
\begin{equation*}
f=\sum_{\alpha=(i_1, \ldots ,i_n)}c_{\alpha}x_1^{i_1}\ldots
x_n^{i_n}.
\end{equation*}
we formally associate the differential operator $D_f$:
\begin{equation*}
D_f=\sum_{\alpha=(i_1, \ldots
,i_n)}c_{\alpha}\frac{\partial^{i_1}}{\partial x_1^{i_1}}\cdots
\frac{\partial^{i_n}}{\partial x_n^{i_n}}.
\end{equation*}
We define the following metric on $P_{n,2k}$, which we call the
differential metric:
\begin{equation*}
\inprodd{f}{g}=D_f(g).
\end{equation*}
It is not hard to check that this indeed defines a symmetric
positive definite bilinear form, which is invariant under the
action of $SO(n)$. For a point $v \in \sph$ we will use $v^{2k}$
to denote the polynomial
\begin{equation*}
v^{2k}=(v_1x_1+ \ldots +v_nx_n)^{2k}.
\end{equation*}
\indent We also define an important linear operator $T: P_{n,2k}
\to P_{n,2k}$, which to a form $f \in P_{n,2k}$ associates
weighted average of forms $v^{2k}$ with the weight $f(v)$:
\begin{equation*}
T(f)=\int_{\sph} f(v)v^{2k} \, d\sigma(v).
\end{equation*}
The operator $T$ was first introduced in a very different form by
Reznick in \cite{rez2}; we take our definition from \cite{greg2}.
The operator $T$ acts as a switch between our standard integral
metric and the differential metric in the following sense:
\begin{lemma}
\label{scalarprodchange} The following identity relating the
operator $T$ and the two metrics holds,
\begin{equation*}
\inprodd{Tf}{g}=(2k)!\inprod{f}{g}.
\end{equation*}
\end{lemma}
\begin{proof}
We observe that
\begin{equation*}
\inprodd{Tf}{g}=\inprodd{\int_{\sph} f(v)v^{2k} \,
d\sigma(v)}{g}=\int_{\sph} \inprodd{f(v)v^{2k}}{g} \, d \sigma(v).
\end{equation*}
Since
\begin{equation*}
\inprodd{v^{2k}}{g}=(2k)!g(v),
\end{equation*}
it follows that
\begin{equation*}
\inprodd{Tf}{g}=(2k)!\int_{\sph} f(v)g(v) \, d
\sigma(v)=(2k)!\inprod{f}{g}.
\end{equation*}
\end{proof}
\indent Let $L$ be a full-dimensional cone in $P_{n,2k}$ such that
$r^{2k}$ is in the interior of $L$ and $\int_{\sph} f \, d\sigma >
0$ for all non-zero $f$ in $L$. We define $\widetilde{L}$ as the
set of all forms $f$ in $M$ such that $f+r^{2k}$ lies in $L$,
\begin{equation*}
\widetilde{L}=\{f \in M \mid f+r^{2k} \in L \}.
\end{equation*}
We let $L^*_i$ be the dual cone of $L$ in the integral metric and
$L^*_d$ be the dual cone of $L$ in the differential metric.
\begin{eqnarray*}
L^*_i=\{f \in P_{n,2k} \mid \inprod{f}{g} \geq 0 \foral g \in
L\}, \\
L^*_d=\{f \in P_{n,2k} \mid \inprodd{f}{g} \geq 0 \foral g \in
L\}.
\end{eqnarray*}
We observe that $r^{2k}$ is in the interior of both $L^*_i$ and
$L^*_d$ and also $\int_{\sph} f \, d\sigma > 0$ for all non-zero
$f$ in both of the dual cones. Therefore we can similarly define
$\widetilde{L^*_i}$ and $\widetilde{L^*_d}$ as sets of all forms
$f$ in $M$ such that $f+r^{2k}$ lies in the respective cone.
\begin{lemma}
\label{volumeswitch} Let $L$ be a full-dimensional cone in
$P_{n,2k}$ such that $r^{2k}$ is the interior of $L$ and
$\int_{\sph} f \, d\sigma
> 0$ for all $f$ in $L$. Then there is the following relationship
between the volumes of $\widetilde{L^*_i}$ and $\widetilde{L^*_d}$
\begin{equation*}
\frac{k!}{(n/2+2k)^{k}} \leq \bigg{(}\frac{\text{Vol} \,
\widetilde{L^*_d}}{\text{Vol}\, \widetilde{L^*_i}}
\bigg{)}^{1/D_M} \leq
\bigg{(}\frac{k!}{(n/2+k)^{k}}\bigg{)}^{\alpha},
\end{equation*}
where
\begin{equation*}
\alpha=1-\bigg{(}\frac{2k-1}{2k+n-2}\bigg{)}^2.
\end{equation*}

\end{lemma}
\begin{proof}
From Lemma \ref{scalarprodchange} we see that
\begin{equation*}
\inprod{f}{g} \geq 0 \quad \text{if and only if} \quad
\inprodd{Tf}{g} \geq 0 \quad \text{for all} \quad f,g \in
P_{n,2k}.
\end{equation*}
Therefore it follows that $T$ maps $L_i^*$ to $L_d^*$,
\begin{equation*}
T(L_i^*)=L_d^*.
\end{equation*}
It is hot hard to show that
\begin{equation*}
T(r^{2k})=cr^{2k} \quad \text{where} \quad c=\int_{\sph} x_1^{2k}
\,
d\sigma=\frac{\Gamma(\frac{2k+1}{2})\Gamma(\frac{n}{2})}{\sqrt{\pi}\Gamma(\frac{n+2k}{2})}.
\end{equation*}
Therefore $\frac{1}{c}T$ fixes the hyperplane of all forms of
integral 1 on the sphere and therefore $\frac{1}{c}T$ maps the
section $\widetilde{L_i^*}$ to $\widetilde{L_d^*}$. \\
\indent It is possible to describe precisely the action of
$\frac{1}{c}T$ on $M_{n,2k}$, see \cite{greg2}. It can be shown
that $\frac{1}{c}T$ is a contraction operator and the exact
coefficients of contraction can be computed. We only need the
following estimate, which follows from \cite{greg2} Lemma 7.4 by
estimating the change in volume to be at most the largest
contraction coefficient:
\begin{equation*}
\bigg{(}\frac{\text{Vol}\, \widetilde{L_d^*}}{\text{Vol}\,
\widetilde{L_i^*}}\bigg{)}^{1/D_M} \geq
\frac{k!\Gamma(k+n/2)}{\Gamma(2k+n/2)}.
\end{equation*}
We observe that
\begin{equation*}
\frac{k!\Gamma(k+n/2)}{\Gamma(2k+n/2)} \geq
\frac{k!}{(n/2+2k)^{k}},
\end{equation*}
and therefore,
\begin{equation*}
\bigg{(}\frac{\text{Vol}\, \widetilde{L_d^*}}{\text{Vol}\,
\widetilde{L_i^*}}\bigg{)}^{1/D_M} \geq \frac{k!}{(n/2+2k)^{k}}.
\end{equation*}
\indent Also from Lemma 7.4 of \cite{greg2} it follows that
contraction by the largest coefficient occurs in the space of all
harmonic polynomials of degree $2k$ which has dimension
\begin{equation*}
D_H=\binom{n+2k-1}{2k}-\binom{n+2k-3}{2k-2}.
\end{equation*}
Since the dimension of the ambient space $M$ is
\begin{equation*}
D_M=\binom{n+2k-1}{2k}-1,
\end{equation*}
we can estimate that
\begin{equation*}
\frac{D_H}{D_M} \geq 1-\bigg{(}\frac{2k-1}{n+2k-2}\bigg{)}^2.
\end{equation*}
Since we can also estimate the largest contraction coefficient
from above,
\begin{equation*}
\frac{k!\Gamma(k+n/2)}{\Gamma(2k+n/2)} \leq
\frac{k!}{(n/2+k)^{k}},
\end{equation*}
the theorem now follows.
\end{proof}
\indent We also show the following theorem, which allows us to
compare the cone of sums of squares to its dual.
\begin{lemma}
\label{dualinclusion} The dual cone to the cone of sums of squares
in the differential metric $Sq_d^*$ is contained in the cone of
sums of squares $Sq$,
\begin{equation*}
Sq_d^* \subseteq Sq.
\end{equation*}
\end{lemma}
\begin{proof}
In this proof we will work exclusively with the differential
metric on $P_{n,k}$ and $P_{n,2k}$. Let $W$ be the space of
quadratic forms on $P_{n,k}$. For $A,\,B$ in $W$, with
corresponding symmetric matrices $M_A, \, M_B$ the inner product
of $A$ and $B$ is given by,
\begin{equation*}
\inprod{A}{B}=\text{tr} \, M_AM_B.
\end{equation*}
For $q \in P_{n,k}$ let $A_q$ be the rank one quadratic form
giving the square of the inner product with $q$:
\begin{equation*}
A_q(p)=\inprodd{p}{q}^2.
\end{equation*}
Then for any $B \in W$
\begin{equation*}
\inprod{A_q}{B}=B(q).
\end{equation*}
\indent Now suppose $f \in Sq_d^*$. Let $H_f$ be the following
quadratic form on $P_{n,k}$:
\begin{equation*}
H_f(p)=\inprodd{p}{f^2}.
\end{equation*}
Since $f \in Sq_d^*$, the quadratic form $H_f$ is clearly positive
semidefinite. Therefore $H_f$ can be written as a nonnegative
linear combination of forms of rank 1:
\begin{equation}
\label{decompose} H_f=\sum A_q \qquad \text{for some} \qquad q \in
P_{n,k}.
\end{equation}
\indent Let $V$ be the subspace of $W$ given by the linear span of
the forms $H_f$ for all $f \in P_{n,2k}$. Let $\mathbb{P}$ be the
operator of orthogonal projection onto $V$. We claim that
\begin{equation*}
\mathbb{P}(A_q)=\binom{2k}{k}^{-1}H_{q^2}.
\end{equation*}
It suffices to show that $A_q-\binom{2k}{k}^{-1}H_{q^2}$ is
orthogonal to the forms $H_{v^{2k}}$ since these forms span $V$.
We observe that
\begin{equation*}
H_{v^{2k}}(p)=(2k)!p(v)^{2k}=\frac{(2k)!A_{v^k}(p)}{(k!)^2}=\binom{2k}{k}A_{v^k}(p).
\end{equation*}
Therefore we see that
\begin{equation*}
\inprod{A_q-\binom{2k}{k}^{-1}\! H_{q^2}}{H_{v^{2k}}}= \!
H_{v^{2k}}(q)-
\inprod{H_{q^2}}{A_{v^k}}=H_{v^{2k}}(q)-H_{q^2}(v^k)\!=\!0.
\end{equation*}
\indent Now we apply $\mathbb{P}$ to both sides of
\eqref{decompose}. It follows that
\begin{equation*}
H_f=\mathbb{P}\bigg(\sum A_q \bigg)=\sum
\binom{2k}{k}^{-1}H_{q^2}=\binom{2k}{k}^{-1}H_{\sum q^2}.
\end{equation*}
Therefore $f$ is a sum of squares.
\end{proof}
\section{Sums of Squares}
\setcounter{equation}{0}

In this section we prove Theorem \ref{squarenormbound}. The full
statement of the bounds is the following,
\begin{theorem}
\label{squarenormboundfull} There are the following bounds for the
volume of $\widetilde{Sq}$:
\begin{equation*}
\frac{(k!)^2}{4^{2k}(2k)!\sqrt{24}}\frac{n^{k/2}}{(n/2+2k)^{k}}
\leq \bigg{(}\frac{\text{Vol}\, \widetilde{Sq}}{\text{Vol}\,
B_M}\bigg{)}^{1/D_M} \leq \frac{4^{2k}(2k)!\sqrt{24}}{k!}\,
n^{-k/2}.
\end{equation*}
\end{theorem}

\subsection{Proof of the Upper Bound}\hspace{1mm}\\
\indent Let us begin by considering the support function of
$\widetilde{Sq}$, which we call $L_{\widetilde{Sq}}$:
\begin{equation*}
L_{\widetilde{Sq}}(f)=\max_{g \, \in \, \widetilde{Sq}} \,
\inprod{f}{g}.
\end{equation*}
The average width $W_{\widetilde{Sq}}$ of $\widetilde{Sq}$ is
given by
\begin{equation*}
W_{\widetilde{Sq}}=2\int_{S_M} L_{\widetilde{Sq}} \, d\mu.
\end{equation*}
We now recall Urysohn's Inequality \cite[p.318]{schneider} which
applied to $\widetilde{Sq}$ gives
\begin{equation}
\label{urysohn} \bigg{(}\frac{\text{Vol} \,
\widetilde{Sq}}{\text{Vol} \, B_M}\bigg{)}^{\frac{1}{D_M}} \, \leq
\frac{W_{\widetilde{Sq}}}{2}.
\end{equation}
Therefore it suffices to obtain an upper bound for
$W_{\widetilde{Sq}}$. \\
\indent Let $S_{P_{n,k}}$ denote the unit sphere in $P_{n,k}$. We
observe that extreme points of $\widetilde{Sq}$ have the form
\begin{equation*}
g^2-r^{2k} \qquad \text{where} \qquad g \in P_{n,k} \qquad
\text{and} \qquad \int_{\sph}g^2 \,d\sigma=1.
\end{equation*}
For $f \in M$,
\begin{equation*}
\inprod{f}{r^{2k}}=\int_{\sph}f \,d\sigma=0,
\end{equation*}
and therefore,
\begin{equation*}
L_{\widetilde{Sq}}(f)=\max_{g \, \in S_{P_{n,k}}} \inprod{f}{g^2}.
\end{equation*}
\indent We now introduce a norm on $P_{n,2k}$, which we denote $||
\ ||_{sq}$:
\begin{equation*}
||f||_{sq}=\max_{g \, \in \, S_{P_{n,k}}} |\inprod{f}{g^2}|.
\end{equation*}
It is clear that
\begin{equation*}
L_{Sq}(f) \leq ||f||_{Sq}.
\end{equation*}
Therefore by \eqref{urysohn} it follows that
\begin{equation*}
\bigg{(}\frac{\text{Vol} \, \widetilde{Sq}}{\text{Vol} \,
B_M}\bigg{)}^{\frac{1}{D_M}} \, \leq \int_{S_M} ||f||_{sq} \, d
\mu.
\end{equation*}
The proof of the upper bound of Theorem \ref{squarenormboundfull}
is reduced to the estimate below.
\begin{theorem}
\label{squareballest} There is the following bound for the average
$|| \ ||_{sq}$ over $S_M$:
\begin{equation*}
\int_{S_M} ||f||_{sq} \, d\mu \, \leq \,
\frac{4^{2k}(2k)!\sqrt{24}}{k!}\, n^{-k/2}.
\end{equation*}
\end{theorem}
\begin{proof}
For $f \in P_{n,2k}$ we introduce a quadratic form $H_f$ on
$P_{n,k}$:
\begin{equation*}
H_f(g)=\inprod{f}{g^2} \qquad \text{for} \qquad g \in P_{n,k}.
\end{equation*}
We note that
\begin{equation*}
||f||_{sq}=\max_{g \, \in \,
S_{P_{n,k}}}|\inprod{f}{g}|=||H_f||_{\infty}.
\end{equation*}
We bound $||H_f||_{\infty}$ by a high $L^{2p}$ norm of $H_f$.
Since $H_f$ is a form of degree 2 on the vector space $P_{n,k}$ of
dimension $D_{n,k}$ it follows by the inequality of Barvinok in
\cite{barv} applied in the same way as in the proof of Theorem
\ref{posmain} that
\begin{equation*}
||H_f||_{\infty} \leq 2\sqrt{3} \, ||H_f||_{2D_{n,k}}.
\end{equation*}
Therefore it suffices to estimate:
\begin{equation*}
A= \int_{S_M} ||H_f||_{2D_{n,k}} \, d\mu = \int_{S_M}
\bigg{(}\int_{S_{P_{n,k}}} \inprod{f}{g^2}^{\, 2D_{n,k}} \, d
\sigma(g) \, d\mu(f) \bigg{)}^{\frac{1}{2D_{n,k}}}.
\end{equation*}
We apply H\"{o}lder's inequality to see that
\begin{equation*}
A \leq \bigg{(}\int_{S_M} \int_{S_{P_{n,k}}} \inprod{f}{g^2}^{\,
2D_{n,k}} \, d \sigma(g) \, d\mu(f) \bigg{)}^{\frac{1}{2D_{n,k}}}.
\end{equation*}
By interchanging the order of integration we obtain
\begin{equation}
\label{weird1} A \leq \bigg{(}\int_{S_{P_{n,k}}} \int_{S_M}
\inprod{f}{g^2}^{\, 2D_{n,k}} \, d \mu(f) \, d\sigma(g)
\bigg{)}^{\frac{1}{2D_{n,k}}}.
\end{equation}
\indent Now we observe that the inner integral
\begin{equation*}
\int_{S_M} \inprod{f}{g^2}^{\, 2D_{n,k}} \, d \mu(f),
\end{equation*}
clearly depends only on the length of the projection of $g^2$ into
$M$. Therefore we have
\begin{equation*}
\int_{S_M} \inprod{f}{g^2}^{\, 2D_{n,k}} \, d \mu(f) \leq \,
||g^2||_2^{2D_{n,k}}\int_{S_M}\inprod{f}{p}^{2D_{n,k}} \, d\mu(f),
\end{equation*}
\begin{equation*}
\text{for any} \quad p \in S_M.
\end{equation*}
We observe that
\begin{equation*}
||g^2||_2=(||g||_4)^2 \qquad \text{and} \qquad ||g||_2=1.
\end{equation*}
By a result of Duoandikoetxea \cite{duo} Corollary 3 it follows
that
\begin{equation*}
||g^2||_2 \leq 4^{2k}.
\end{equation*}
Hence we obtain
\begin{equation*}
\int_{S_M} \inprod{f}{g^2}^{\, 2D_{n,k}} \, d \mu(f) \leq
4^{4kD_{n,k}} \int_{S_V}\inprod{f}{p}^{2D_{n,k}} \, d\mu(f).
\end{equation*}
We note that this bound is independent of $g$ and substituting
into \eqref{weird1} we get
\begin{equation*}
A \leq 4^{2k}\bigg{(}\int_{S_V}\inprod{f}{p}^{2D_{n,k}} \,
d\mu(f)\bigg{)}^{\frac{1}{2D_{n,k}}}.
\end{equation*}
\indent Since $p \in S_M$ we have
\begin{equation*}
\int_{S_M}\inprod{f}{p}^{2D_{n,k}} \, d\mu(f) =
\frac{\Gamma(D_{n,k}+\frac{1}{2})\Gamma(\,
\frac{1}{2}D_M)}{\sqrt{\pi} \, \Gamma(D_{n,k}+\frac{1}{2}D_M)}.
\end{equation*}
We use the following easy inequalities:
\begin{equation*}
\bigg{(}\frac{\Gamma(\,
\frac{1}{2}D_M)}{\Gamma(D_{n,k}+\frac{1}{2}D_M)}\bigg{)}^{\frac{1}{2D_{n,k}}}
 \leq \sqrt{\frac{2}{D_M}}
\end{equation*}
and
\begin{equation*}
\bigg{(}\frac{\Gamma(D_{n,k}+\frac{1}{2})}{\sqrt{\pi}}\bigg{)}^{\frac{1}{2D_{n,k}}}
\leq \sqrt{D_{n,k}},
\end{equation*}
to see that
\begin{equation*}
A \leq 4^{2k}\sqrt{\frac{2D_{n,k}}{D_M}}.
\end{equation*}
We now recall that
\begin{equation*}
D_{n,k}=\binom{n+k-1}{k} \qquad \text{and} \qquad
D_M=\binom{n+2k-1}{2k}-1.
\end{equation*}
Therefore
\begin{equation*}
\sqrt{\frac{D_{n,k}}{D_M}} \, \leq  \, \frac{(2k)!}{k!} \,
n^{-k/2}.
\end{equation*}
Thus
\begin{equation*}
A \leq \frac{4^k(2k)!\sqrt{2}}{k!} \, n^{-k/2}.
\end{equation*}
The theorem now follows.
\end{proof}
\subsection{Proof of the Lower Bound}\hspace{1mm}\\
\indent We begin with a corollary of Theorem \ref{squareballest}.
Let $B_{sq}$ be the unit ball of the norm $|| \ ||_{sq}$,
\begin{equation*}
B_{sq}=\{ f \in M \mid ||f||_{sq} \leq 1 \}.
\end{equation*}
From Theorem \ref{squareballest} we know that
\begin{equation*}
\int_{S_M} ||f||_{sq} \, d\mu \, \leq \,
\frac{4^{2k}(2k)!\sqrt{24}}{k!}\, n^{-k/2}.
\end{equation*}
It follows in the same way as in the section 3.1 that
\begin{equation*}
\bigg{(}\frac{\text{Vol}\, B_{sq}}{\text{Vol}\,
B_M}\bigg{)}^{1/D_M} \geq \frac{k!}{4^{2k}(2k)!\sqrt{24}}\,
n^{k/2}.
\end{equation*}
\indent Now let $\widetilde{Sq}^{\circ}$ be the polar of
$\widetilde{Sq}$ in $M$. It follows easily that $B_{sq}$ is the
intersection of $\widetilde{Sq}^{\circ}$ and
$-\widetilde{Sq}^{\circ}$.
\begin{equation*}
B_{sq}=\widetilde{Sq}^{\circ} \cap -\widetilde{Sq}^{\circ}.
\end{equation*}
Let ${Sq^*_i}$ be the dual cone of $Sq$ in the integral metric and
let $\widetilde{Sq^*_i}$ be defined in the same way as for the
previous cones. It is not hard to check that
$\widetilde{Sq}^{\circ}$ is the negative of $\widetilde{Sq^*_i}$,
\begin{equation*}
\widetilde{Sq}^{\circ}=-\widetilde{Sq^*_i}.
\end{equation*}
Therefore we see that
\begin{equation*}
\label{short1} \bigg{(}\frac{\text{Vol}\,
\widetilde{Sq^*_i}}{\text{Vol}\, B_M}\bigg{)}^{1/D_M} \geq
\frac{k!}{4^{2k}(2k)!\sqrt{24}}\, n^{k/2}.
\end{equation*}
Now we observe that $r^{2k}$ is in the interior of $Sq$ and also
for all non-zero $f$ in $Sq$ we have $\int_{\sph} f \, d\sigma >
0$. Therefore we can apply Lemma \ref{volumeswitch} to $Sq$ and it
follows that
\begin{equation*}
\bigg{(}\frac{\text{Vol}\, \widetilde{Sq^*_d}}{\text{Vol}\,
\widetilde{Sq^*_i}}\bigg{)}^{1/D_M} \geq \frac{k!}{(n/2+2k)^{k}}.
\end{equation*}
Combining with \eqref{short1} we see that
\begin{equation*}
\bigg{(}\frac{\text{Vol}\, \widetilde{Sq^*_d}}{\text{Vol}\,
B_M}\bigg{)}^{1/D_M} \geq
\frac{(k!)^2}{4^{2k}(2k)!\sqrt{24}}\frac{n^{k/2}}{(n/2+2k)^{k}}.
\end{equation*}
By Lemma \ref{dualinclusion} we know that $Sq_d^*$ in contained in
$Sq$ and therefore
\begin{equation*}
\widetilde{Sq_d^*} \subseteq \widetilde{Sq}.
\end{equation*}
The lower bound now follows.

\section{Sums of 2k-th Powers of Linear Forms}
\setcounter{equation}{0}

In this section we prove Theorem \ref{powersnormbound}. Here is
the precise statement of the bounds,
\begin{theorem}
\label{powersnormboundfull} There are the following bounds for the
volume of $\widetilde{Lf}$:
\begin{equation*}
\frac{k!\sqrt{4k^2+n-2}}{4k\sqrt{2}(n/2+2k)^k} \leq \!
\bigg{(}\frac{\text{Vol} \, \widetilde{Lf}}{\text{Vol} \,
B_M}\bigg{)}^{1/D_M} \! \!  \leq
2\sqrt{n(4k+2)}\bigg{(}\frac{k!}{(n/2+k)^{k}}\bigg{)}^{\alpha},
\end{equation*}
where
\begin{equation*}
\alpha=1-\bigg{(}\frac{2k-1}{n+2k-2}\bigg{)}^2.
\end{equation*}
\end{theorem}

\subsection{Proof of the Lower Bound} \hspace{1mm}\\
\indent We observe that the cone of sums of $2k$-th powers of
linear forms is dual to the cone of nonnegative polynomials in the
differential metric,
\begin{equation*}
Lf=C^*_d,
\end{equation*}
since in the differential metric,
\begin{equation*}
\inprodd{f}{v^{2k}}=(2k)!f(v) \foral f \in P_{n,2k}.
\end{equation*}
Therefore it follows that
\begin{equation*}
\widetilde{Lf}=\widetilde{C^*_d}.
\end{equation*}
\indent We first consider the dual cone $C^*_i$ of $C$ in the
integral metric. Similarly to the situation with the cone of sums
of squares it is not hard to check that the dual
$\widetilde{C}^{\circ}$ of $\widetilde{C}$ in $M$ with respect to
the integral metric is $-\widetilde{C^*_i}$,
\begin{equation*}
\widetilde{C}^{\circ}=-\widetilde{C^*_i}.
\end{equation*}
We recall that in Section 3.2 we have shown \eqref{ratio}:
\begin{equation*}
\bigg{(}\frac{\text{Vol} \, \widetilde{C}^{\circ}}{\text{Vol} \,
B_M}\bigg{)}^{1/D_M} \geq
\frac{1}{4}\bigg{(}\frac{4k^2+n-2}{2k^2}\bigg{)}^{1/2}.
\end{equation*}
Since $C$ has $r^{2k}$ in its interior and $\int_{\sph} f \,
d\sigma > 0$ for all non-zero $f$ in $C$, we can apply Lemma
\ref{volumeswitch} to $C$ and we obtain,
\begin{equation*}
\bigg{(}\frac{\text{Vol} \, \widetilde{C^*_d}}{\text{Vol}\,
\widetilde{C^*_i}} \bigg{)}^{1/D_M} \geq \frac{k!}{(n/2+2k)^{k}}.
\end{equation*}
Since $\widetilde{Lf}=\widetilde{C^*_d}$ and
$\widetilde{C}^{\circ}=-\widetilde{C^*_i}$ we can combine with
\eqref{ratio} and we get:
\begin{equation*}
\bigg{(}\frac{\text{Vol} \, \widetilde{Lf}}{\text{Vol} \,
B_M}\bigg{)}^{1/D_M} \geq
\frac{k!}{4k\sqrt{2}}\frac{(4k^2+n-2)^{1/2}}{(n/2+2k)^k}.
\end{equation*}

\subsection{Proof of the Upper Bound} \hspace{1mm} \\
\indent We begin by applying the Blaschke-Santal\'{o} inequality
to $\widetilde{C}$ as in Section 3.2 to obtain
\begin{equation*}
\frac{\text{Vol}\, \widetilde{C} \, \text{Vol}\,
\widetilde{C}^{\circ}}{(\text{Vol}\, B_M)^2} \leq 1.
\end{equation*}
Since $\widetilde{C}^{\circ}=-\widetilde{C^*_i}$ we can rewrite
this to get
\begin{equation*}
\bigg{(}\frac{\text{Vol}\, \widetilde{C^*_i}}{\text{Vol}\,
B_M}\bigg{)}^{1/D_M} \leq \bigg{(}\frac{\text{Vol}\,
B_M}{\text{Vol}\, \widetilde{C}}\bigg{)}^{1/D_M}.
\end{equation*}
We observe that by the lower bound of Theorem \ref{posmainfull} it
follows that
\begin{equation}
\label{short2} \bigg{(}\frac{\text{Vol}\,
\widetilde{C^*_i}}{\text{Vol}\, B_M}\bigg{)}^{1/D_M} \leq
2\sqrt{n(4k+2)}.
\end{equation}
Now we apply the upper bound of Lemma \ref{volumeswitch} to $C$
and we get
\begin{equation*}
\bigg{(}\frac{\text{Vol} \, \widetilde{C^*_d}}{\text{Vol}\,
\widetilde{C^*_i}} \bigg{)}^{1/D_M} \leq
\bigg{(}\frac{k!}{(n/2+k)^{k}}\bigg{)}^{\alpha},
\end{equation*}
where
\begin{equation*}
\alpha=1-\bigg{(}\frac{2k-1}{n+2k-2}\bigg{)}^2.
\end{equation*}
The upper bound now follows by combining with \eqref{short2}.

\medskip
\textsc{Department of Mathematics, University of Michigan, \\Ann
Arbor, MI 48109-1109, USA}\\
\textit{Email address:} gblekher@umich.edu

\begin{thebibliography}{dmst}
\bibitem[1]{barv} A.I. Barvinok,
\textit{Estimating $L^{\infty}$ norms by $L^{2k}$ norms for
functions on orbits}. Foundations of Computational Mathematics, 2
(2002),  no. 4, 393-412.

\bibitem[2]{greg2} G. Blekherman
\textit{Convexity properties of the cone of nonnegative polynomials}, arXiv preprint math.CO/0211176 (2002),
Discrete and Computational Geometry to appear.


\bibitem[3]{greg} G. Blekherman \textit{There are significantly
more nonnegative polynomials than sums of squares}, arXiv preprint
math.AG/0309130 (2003).

\bibitem[4]{choi} M. D. Choi,  T. Y. Lam, B. Reznick,
\textit{Even symmetric sextics.}  Math. Z.  195  (1987),  no. 4,
559-580.

\bibitem[5]{duo} J. Duoandikoetxea,
\textit{Reverse H\"{o}lder inequalities for spherical harmonics.}
Proc. Amer. Math. Soc. 101 (1987), no. 3, 487-491.

\bibitem[6]{hardy} G. H. Hardy, J. E. Littlewood, G. P\'{o}lya,
\textit{Inequalities.} Reprint of the 1952 edition. Cambridge
Mathematical Library. Cambridge University Press, Cambridge, 1988.

\bibitem[7]{hilbert} D. Hilbert,
\textit{\"{U}ber die Darstellung definiter Formen als Summe von
Formenquadraten.} Math. Ann. 32, 342-350 (1888). Ges Abh. vol. 2,
415-436. Chelsea Publishing Co., New York, (1965).


\bibitem[8]{kellogg} O. Kellogg, \textit{On bounded polynomials in several
variables.} Math. Z. 27, 1928, 55-64.


\bibitem[9]{santalo} M. Meyer, A. Pajor.
\textit{On
 the Blaschke-Santaló inequality.}  Arch. Math. (Basel)
55 (1990),  no. 1, 82-93.


\bibitem[10]{pach} J. Pach, P. Agarwal. \textit{Combinatorial
Geometry.} Wiley-Interscience Series in Discrete Mathematics and
Optimization. John Wiley \& Sons, Inc., New York, 1995.


\bibitem[11]{pisier}G. Pisier, \textit{The Volume of Convex Bodies and Banach space
Geometry.} Cambridge Tracts in Mathematics, 94. Cambridge
University Press, Cambridge, 1989.

\bibitem[12]{rez} B. Reznick,
\textit{Sums of even powers of real linear forms}, Mem. Amer.
Math. Soc. 96 (1992), no. 463.

\bibitem[13]{rez2} B. Reznick, \textit{Uniform denominators in Hilbert's
seventeenth problem.} Math. Zeitschrift. 220 (1995), no. 1,
75--97.


\bibitem[14]{rez4}B. Reznick,
\textit{Some concrete aspects of Hilbert's 17th Problem.} Contemp.
Math., 253 (2000), 251-272.

\bibitem[15]{schneider} R. Schneider, \textit{Convex bodies: the Brunn-Minkowski
theory.} Encyclopedia of Mathematics and its Applications, 44.
Cambridge University Press, Cambridge, 1993.

\bibitem[16]{parrilo}P. A. Parrilo, B. Sturmfels. \textit{Minimizing polynomials
functions.} Submitted to the DIMACS volume of the Workshop on
Algorithmic and Quantitative Aspects of Real Algebraic Geometry in
Mathematics and Computer Science.


\bibitem[17]{vilenkin} N. Ja. Vilenkin,
\textit{Special Functions and the Theory of Group
Representations.} Translations of Mathematical Monographs, Vol.
22, American Mathematical Society (1968).
\end{thebibliography}
\end{document}